\documentclass[12pt]{amsart}
\usepackage{amsfonts, amsmath, latexsym, epsfig}
\usepackage{amssymb}
\usepackage{epsf}
\usepackage{color}
\usepackage{url}

\newcommand{\RR}{\ensuremath{\mathbb{R}}}
\newcommand{\NN}{\ensuremath{\mathbb{N}}}

\newcommand{\ZZ}{\ensuremath{\mathbb{Z}}}

\newtheorem{proposition}{Proposition}
\newtheorem{theorem}{Theorem}

\newtheorem{lemma}{Lemma}

\newtheorem{conjecture}{Conjecture}

\newtheorem{definition}{Definition}
\def\QuotS#1#2{\leavevmode\kern-.0em\raise.2ex\hbox{$#1$}\kern-.1em/\kern-.1em\lower.25ex\hbox{$#2$}}

\newcommand{\SymbVN}{P}

\DeclareMathOperator{\Aut}{Aut}
\DeclareMathOperator{\Sym}{Sym}

\DeclareMathOperator{\GL}{GL}

\DeclareMathOperator{\Stab}{Stab}

\DeclareMathOperator{\cov}{cov}
\DeclareMathOperator{\tvol}{vol}

\begin{document}

\author{Michel Deza}
\address{Michel Deza, \'Ecole Normale Sup\'erieure, Paris} 
\email{Michel.Deza@ens.fr}

\author{Mathieu Dutour Sikiri\'c}
\address{Mathieu Dutour Sikiri\'c, Rudjer Bo\u skovi\'c Institute, Bijeni\u cka 54, 10000 Zagreb, Croatia, Fax: +385-1-468-0245}
\email{mdsikir@irb.hr}

\thanks{Second author has been supported by the Croatian Ministry of Science, Education and Sport under contract 098-0982705-2707.}

\title{Voronoi Polytopes for Polyhedral Norms on Lattices}
\date{}

\maketitle

\begin{abstract}
A {\em polyhedral norm} is a norm $N$ on $\RR^n$ for which the set $N(x)\leq 1$ is a polytope. This covers the case of the $L^1$ and $L^{\infty}$ norms.
We consider here effective algorithms for determining the Voronoi polytope for such norms with a point set being a lattice.
The algorithms, that we propose, use the symmetries effectively in order to compute a decomposition of the space into convex polytopes named {\em $VN$-spaces}.
The Voronoi polytopes and other geometrical information are easily obtained from it.
\end{abstract}

\section{Introduction}\label{SEC_Intro}

The classical theory of Voronoi polytopes takes its roots in the geometry of 
lattices.
For any lattice $L\subset \RR^n$, the Voronoi polytope of a point $v\in L$ is 
defined as
\begin{equation*}
P_V(v)= \left\{ x\in \RR^n \mbox{~:~} \Vert x - w \Vert \geq \Vert x-v\Vert \mbox{~for~all~}w\in L - \{v\} \right\}.
\end{equation*}
The norm $\Vert\cdot \Vert$ is the standard Euclidean norm on $\RR^n$ and the polytope 
$P_V(v)$ is a convex polytope.
A {\em Delaunay polytope} is the convex hull of the set of points closest to a vertex $v$ of the Voronoi polytope. The set of all Delaunay polytopes defines a tiling of the Euclidean space $\RR^n$. There is a one-to-one correspondence between $k$-dimensional faces of Voronoi polytopes and $(n-k)$-dimensional faces of the Delaunay tiling.
Both Voronoi polytopes and Delaunay polytopes define face-to-face tilings of $\RR^n$ and they are are both useful in a number of discrete-geometric questions \cite{BookSchurmann,BookGruber,BookSugihara}.
In \cite{DSV_voronoi} the second author used Delaunay polytopes and lattice symmetries for efficiently computing the Voronoi polytope of many highly symmetric lattices.

For non-Euclidean norms, much of the theory collapses and one has
to adapt to the case considered.
The Voronoi polytopes are no longer convex, but they remain connected.
Here, we consider how one can compute Voronoi polytope for a special kind
of norms that cover both the classical $L^1$ and $L^{\infty}$ norms.
The basic idea is to decompose the Voronoi polytopes into a number of
convex polytopes on which the considered norm behaves nicely.

By a {\em polyhedral norm} on $\RR^n$, we mean a function $N$ of the form
\begin{equation*}
N(x) = N_{{\mathcal L}}(x) = \max_{\ell\in {\mathcal L}} \ell(x).
\end{equation*}
with ${\mathcal L}$ being a finite set of linear forms on $\RR^n$ that spans
the dual of $\RR^n$. Here we have to assume in addition that the linear
forms in ${\mathcal L}$ are rational, i.e. have rational values on rational
vectors.
We assume further that the set ${\mathcal L}$ used to define $N_{{\mathcal L}}$ is minimal.
The function $N$ is called a {\em norm} since it satisfies the following properties:
\begin{enumerate}
\item Triangle inequality: $N(x+y)\leq N(x) + N(y)$ for $x,y\in\RR^n$.
\item $N(x)=0$ is equivalent to $x=0$
\item Positive linearity: for all $x\in\RR^n$ and $\lambda>0$, it holds $N(\lambda x)= \lambda N(x)$.
\end{enumerate}
The norm $N$ will be {\em symmetric}, i.e. $N(x)=N(-x)$, if and only
if ${\mathcal L} = - {\mathcal L}$.
But we do not assume a priori that $N$ is symmetric.
We will use below terms {\em polyhedral norm} ${\mathcal L}$ and 
{\em polyhedral metric} ${\mathcal L}$ for the polyhedral norm, 
generated by ${\mathcal L}$ and corresponding Minkowski metric.
A norm $N$ is polyhedral if and only if the set $P_N=\left\{x\mbox{~:~}N(x)\leq 1\right\}$ is a polytope and it is symmetric if and only if $P_N$ is centrally-symmetric.

The $L^{\infty}$ norm is obtained by taking
\begin{equation}\label{L_Norm_Linfini}
{\mathcal L} = \left\{\pm e_i^* \mbox{~for~} 1\leq i\leq n\right\}.
\end{equation}
with $e_i^*$ being defined by $e_i^*(x_1, \dots, x_n)=x_i$.
The $L^1$ norm is obtained by taking
\begin{equation}\label{L_Norm_L1}
{\mathcal L} = \left\{\sum_{i=1}^n \epsilon_i e_i^* \mbox{~with~} \epsilon_i=\pm 1\right\}.
\end{equation}
A polyhedral norm is called {\em simplicial} if the set of linear forms in ${\mathcal L}$ is linearly independent, i.e. $P_N$ is a simplex.

We want to use the symmetries preserving the lattice $L$ and the norm $N_{\mathcal L}$.
This would allow to compute covering radius and other geometric data.
One highlight of our method is that we are not limited to the two-dimensional
case.

The study of the complexity of the generalized Voronoi algorithm has been proposed in \cite{AgarwalSharir}.
Several complexity results for $L^1$, $L^{\infty}$ and simplicial metrics are obtained in \cite{BoissonatSharir} under the assumption that the point set, for which the Voronoi domains is computed, are in general position. We are not aware of any implementation of those algorithms.
Several algorithms for computing Voronoi diagram on the plane were proposed in \cite{Fu,Jeong,Lee}.
For the three-dimensional case, a randomized algorithm is proposed in \cite{Le} and an algorithm with almost optimal complexity is proposed in \cite{Koltun}.
In \cite{Chew} an algorithm for computing the Voronoi diagram defined by lines is given.
A completely general approximation algorithm is proposed in \cite{Reem_GeneralNormed}. The algorithm is essentially a Monte-Carlo method obtained by tracing rays from each element of the point set.
In \cite{BregmanVoronoi} efficient algorithms are build for Voronoi diagrams obtained from Bregman distance functions and in \cite{Boissonat_Wormser} several algorithms are given for some special distance functions.
On the other hand, polyhedral functions were used in theoretical studies such as \cite{Manjunath}.
Algorithms for Euclidean metrics are too numerous to list.

The general problem considered here is the computation of the Voronoi
polytope for a point set being $\ZZ^n$. With minimal modifications, we
could treat a point set of the form $\cup_{i=1}^m (c_i + \ZZ^n)$ with $c_i\in [0,1[^n$, i.e. crystallographic structures.

In Section \ref{SEC_Geometry} we define the Voronoi polytopes used in this work and explain their geometry.
In Section \ref{SEC_Algorithms}, the enumeration algorithms are developed in details. Those relies on the {\tt zsolve} program for integer enumeration \cite{4ti2}, the {\tt cdd} polyhedral program \cite{cdd} and the implementation is available from \cite{polyhedral}.
In Section \ref{SEC_Applications}, the implementation is applied on the case of the root lattices $\mathsf{A}_n$ and $\mathsf{D}_n$ for the $L^1$ and $L^{\infty}$ norms.

Our algorithm relies in a key way on the hypothesis of rationality of the set ${\mathcal L}$. It is possible that one could dispense from this hypothesis and build efficient and general algorithms.
The next open question would then be the building of parameter space for the Voronoi polytope of polyhedral metric; the first interesting case would be simplicial norms.

\section{Geometry}\label{SEC_Geometry}

In the Figure below, for the sake of clarity, we often choose to replace $\ZZ^n$ by a lattice $L$ and to leave the polyhedral norm invariant. By a change of basis, those lattice changes can be interpreted as polyhedral norm changes.

\begin{figure}
\begin{minipage}[b]{6.0cm}
\centering
\resizebox{40mm}{!}{\rotatebox{90}{\includegraphics{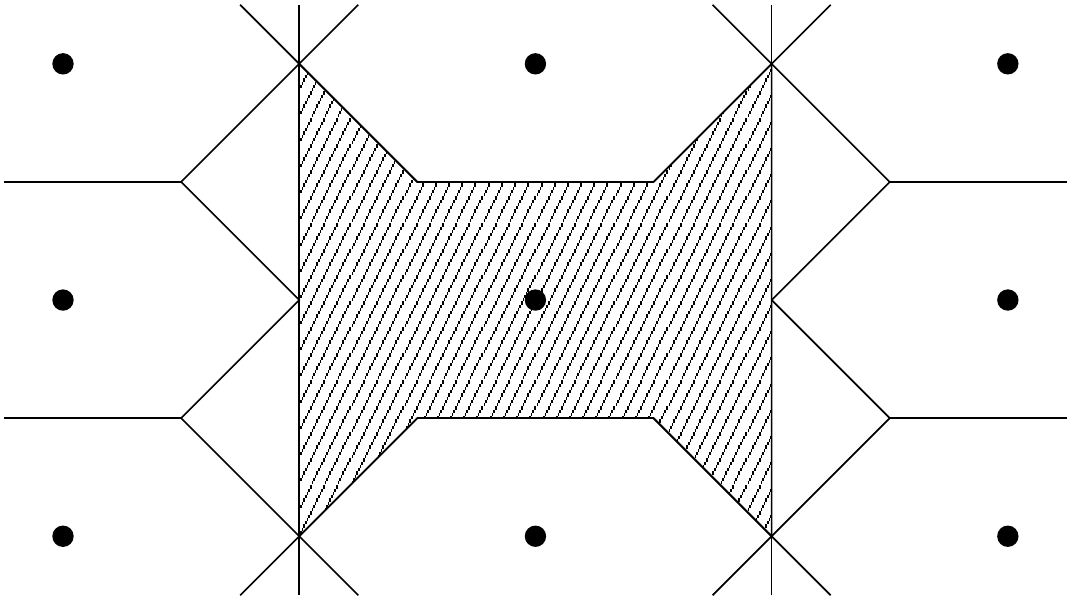}}}\par
a $V_{\leq}$ Voronoi polytope
\end{minipage}
\begin{minipage}[b]{6.0cm}
\centering
\resizebox{40mm}{!}{\rotatebox{90}{\includegraphics{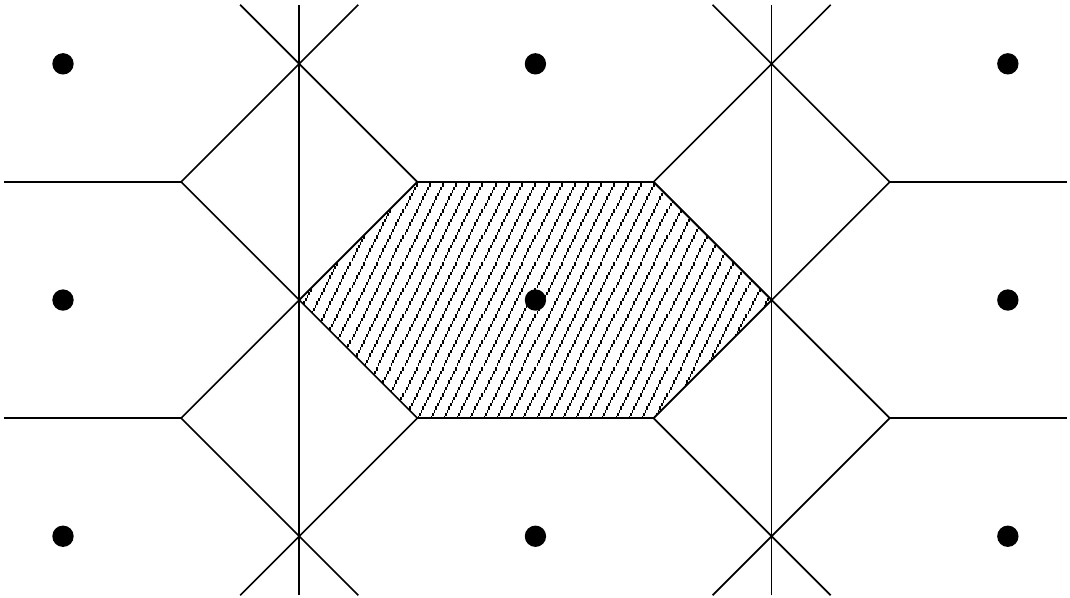}}}\par
a $V_{<}$ Voronoi polytope
\end{minipage}
\begin{center}
\end{center}
\caption{One $V_{\leq}$ and one $V_{<}$ Voronoi polytopes for the $L^{\infty}$ norm.}
\label{VariousVoronoiPolytope}
\end{figure}

Let us take a polyhedral norm $N$. We can define two Voronoi diagrams on the point set $\ZZ^n$:
\begin{equation*}
\begin{array}{rcl}
V_{\leq}(N,v) &=& \left\{ x\in \RR^n \mbox{~:~} N(x-v) \leq N(x-w)\mbox{~for~}w\in\ZZ^n - \{v\}  \right\},\\
V_{ < }(N,v) &=& \left\{ x\in \RR^n \mbox{~:~} N(x-v) < N(x-w)\mbox{~for~}w\in\ZZ^n - \{v\}  \right\}.
\end{array}
\end{equation*}
In the Euclidean case, one usually studies $V_{\leq}$ since
$V_{\leq}$ is the closure of $V_{<}$.
Figure \ref{VariousVoronoiPolytope} shows
that this property does not hold in general.

For the Euclidean metric, the Voronoi polytope is convex.
But only following weaker property holds for polyhedral metrics.

\begin{theorem}\cite{Reem_GeneralNormed}
For any norm $N$ on $\RR^n$, the Voronoi polytopes $V_{\leq}(N,v)$ and $V_{<}(N,v)$ are star convex with respect to $v$.
\end{theorem}
\proof Let us assume $v=0$ and take $x\in V_{\leq}( N,0)$.
If we take $\lambda\in [0,1]$, then we have
\begin{equation*}
\begin{array}{rcl}
N(\lambda x - w) &=& N( -(1-\lambda) x + x-w)\\
&\geq & N(x-w) - N( (1-\lambda) x) \\
&\geq & N(x) - (1-\lambda) N(x) = N(\lambda x)
\end{array}
\end{equation*}
So, $\lambda x\in V_{\leq}(N, 0)$. The proof is similar for $V_{<}(N,v)$. \qed

In \cite{VoronoiII}, a general theory of Voronoi polytopes for Euclidean metrics is developed. As a consequence of this theory, one obtains that as one modifies the Euclidean metric, the Voronoi polytope evolves smoothly.
This property is generalized in \cite{Reem_stability} where a general stability result is proved for {\em uniformly convex spaces} (also called {\em rotund} or {\em strictly convex}), i.e. ones for which the equality $N(x)=N(y)=1$ and $x\not=y$ imply $N(x+y) < 2$. No space, whose norm is defined by a polyhedral metric, is uniformly convex.

No such stability exists, in general, for polyhedral metrics, but one has the following result:

\begin{proposition}
Let ${\mathcal L}_n$ be a sequence of polyhedral metrics that converges towards a polyhedral metric ${\mathcal L}$.
Then one has the following inclusions on the Voronoi polytopes:
\begin{equation*}
V_{<}(N_{{\mathcal L}},v) \subset \underline{\lim}_n V_{<}(N_{{\mathcal L}_n},v)
\mbox{~and~}
\overline{\lim}_n V_{\leq }(N_{{\mathcal L}_n},v) \subset V_{\leq }(N_{{\mathcal L}},v).
\end{equation*}

\end{proposition}
\proof If $x\in V_{<}(N_{{\mathcal L}},v)$, then one has 
$N_{{\mathcal L}}(x-v) < N_{{\mathcal L}}(x-w)$ for $w\in\ZZ^n - \{v\}$.
The Voronoi polytope is bounded; so, only a finite set of those inequalities
is relevant.
For any fixed vector $x$, it holds $\lim_{n\to\infty} N_{{\mathcal L}_n}(x) = N_{{\mathcal L}}(x)$.
As a consequence, for $n$ large enough the inclusion holds.
A similar proof works for the other inclusion. \qed

The above shows that we need to consider both $V_{\leq}(N,v)$ and $V_{ < }(N,v)$ in our work, especially in degenerate situations as defined below:
\begin{definition}
We say that a norm $N$ for $\ZZ^n$ is {\em non-degenerate} if $V_{\leq}(N, v)$ is the closure of $V_{ < }(N, v)$.
\end{definition}

But it is hard to work with the Voronoi polytope directly, and we need instead an object that is more amenable to polyhedral methods.

For a point $x\in \RR^n$, we define the {\em distance to nearest neighbor} as
\begin{equation*}
d_{min}(x,{\mathcal L}) = \min_{v\in \ZZ^n} N_{{\mathcal L}}(x - v).
\end{equation*}
The {\em covering radius} is defined as
\begin{equation*}
\cov({\mathcal L})=\max_{x\in \RR^n} d_{min}(x,{\mathcal L}).
\end{equation*}
In the case of the Euclidean norm $N$, for any two distinct points $v, v'$, the 
set of {\em equidistant points} $x$, i.e., those for which $N(x-v)=N(x-v')$,
is an hyperplane.

\begin{figure}
\begin{minipage}[b]{6.0cm}
\centering
\epsfig{height=3.0cm, file=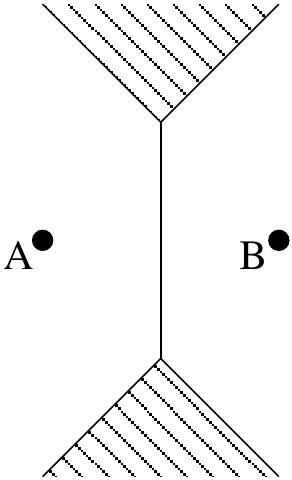}\par
a degenerate configuration
\end{minipage}
\begin{minipage}[b]{6.0cm}
\centering
\epsfig{height=3.0cm, file=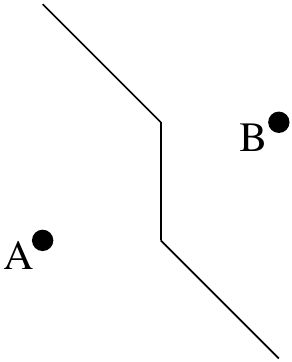}\par
a non-degenerate configuration
\end{minipage}
\begin{center}
\end{center}
\caption{The set of equidistant points for two points and the $L^{\infty}$ norm.}
\label{TwoPoint_Deg_NonDeg}
\end{figure}

\begin{figure}
\begin{center}
\begin{minipage}[b]{6.0cm}
\centering
\epsfig{height=3.5cm, file=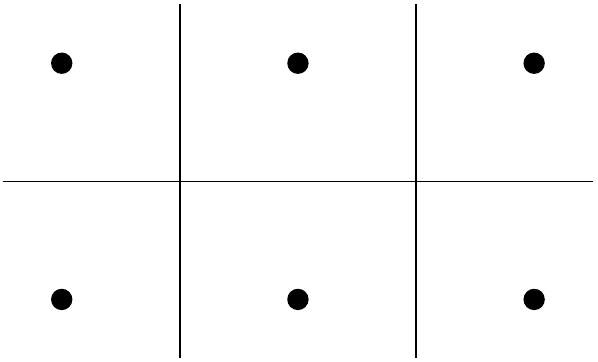}\par
a) a non-degenerate case
\end{minipage}
\begin{minipage}[b]{6.0cm}
\centering
\epsfig{height=3.5cm, file=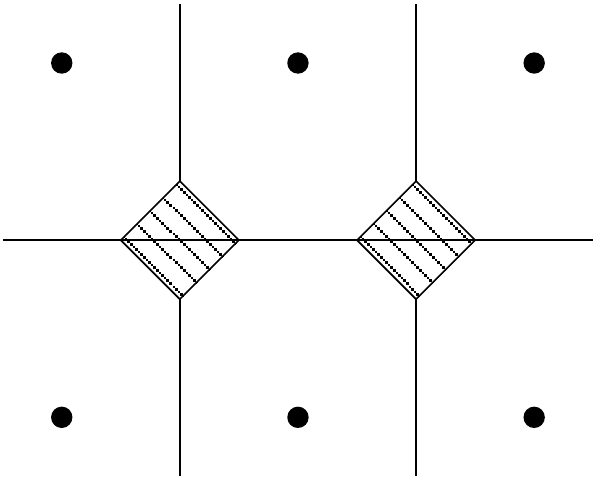}\par
b) a degenerate case
\end{minipage}
\begin{minipage}[b]{6.0cm}
\centering
\epsfig{height=3.5cm, file=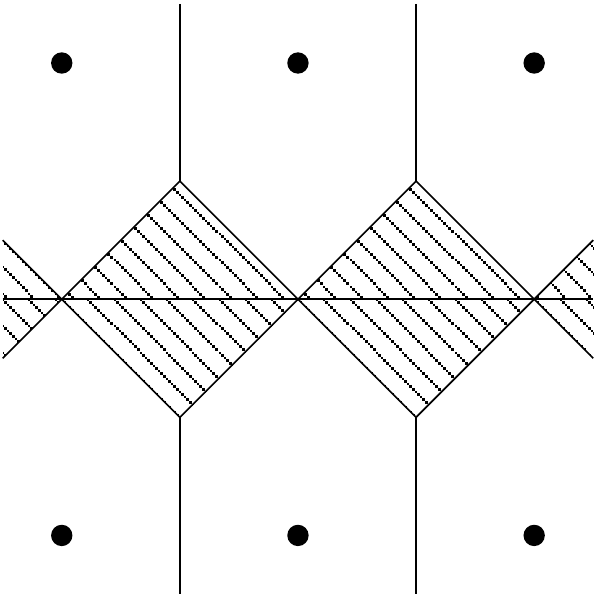}\par
c) a degenerate case\\
\textcolor{white}{Bonjour}
\end{minipage}
\begin{minipage}[b]{6.0cm}
\centering
\epsfig{height=3.5cm, file=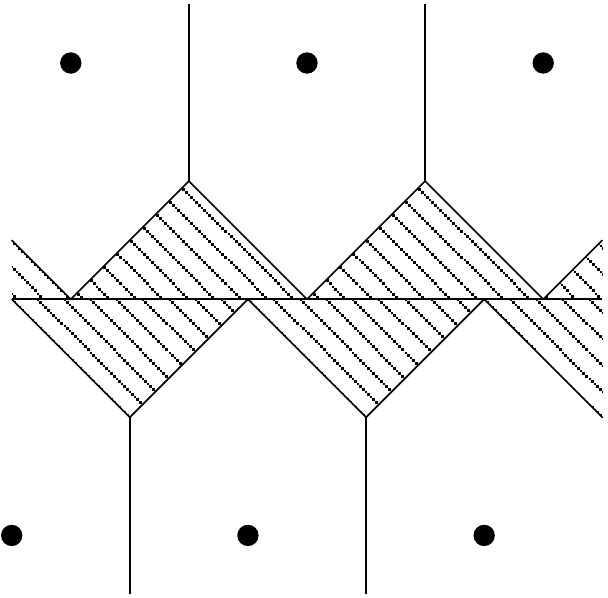}\par
d) a non face-to-face degenerate case
\end{minipage}
\end{center}
\caption{Examples of $VN$-spaces for lattices and the $L^{\infty}$ norm.
Filled areas are the difference $V_{\leq}(N,0) \diagdown V_{<}(N,0)$.}
\label{Lattice_Deg_NonDeg}
\end{figure}

This is no longer true for polyhedral norms. One example is shown on Figure \ref{TwoPoint_Deg_NonDeg}.
For the $L^{\infty}$ norm, we consider the two vectors $A=(-1, 0)$ and $B=(1,0)$.
The points $(x,y)$ with $|y| > 1$ and $|x|\leq |y|-1$ are all at equal distances from $A$ and $B$ and so, the corresponding equidistant points are part of a full-dimensional region.
On the other hand, if we take $\epsilon > 0$ and the points $A=(-1,-\epsilon)$, 
$B=(1,\epsilon)$, then the set of equidistant points are part of an union of segments of dimension $1$.
This phenomenon also occurs when the considered points belong to a lattice (see Figure \ref{Lattice_Deg_NonDeg}).

We now define the main notion of {\em $VN$-space}:
\begin{definition}\label{DEF_VNspace}
Given a polyhedral metric ${\mathcal L}$, an {\em $VN$-space} $P$
is a full-dimensional polytope for which there exist $v\in \ZZ^n$
and $\ell\in {\mathcal L}$ such that:

(i) For all $\ell'\in {\mathcal L}$, the inequality $\ell'(x-v) \leq \ell(x-v)$ is valid for $x\in P$.

(ii) For all $v'\in \ZZ^n$, the inequality $N_{{\mathcal L}}(x-v') \geq \ell(x-v)$ if valid for $x\in P$.

We write $\alpha_P(x)=\ell(x-v)$ and
\begin{equation*}
Near(P)=\{v'\in \ZZ^n\mbox{~:~} N_{{\mathcal L}}(x-v') =\alpha_P(x) \mbox{~for~}x\in P\}
\end{equation*}
\end{definition}

There is a degree of arbitrariness in the above definition of $VN$-space. There is no such thing as canonical $VN$-space associated to a polyhedral metric ${\mathcal L}$.
For example, if we split an $VN$-space $P$ by an hyperplane into two polytopes, then the resulting polytopes are also $VN$-spaces.

Our objective is to tile the space $\RR^n$ with $VN$-spaces, which will allow us to resolve a number of geometrical questions.
Among all such possible $VN$-space decompositions, we are interested in the simplest ones, which will allow easier computations.

The following result follows directly from Definition \ref{DEF_VNspace}:
\begin{proposition}
Given a polyhedral norm ${\mathcal L}$, suppose that we have a
tiling by $VN$-spaces $(P_i)_{i\in I}$.
Then
\begin{equation*}
\cov({\mathcal L})=\max_{i} \max_{x\in P_i} \alpha_{P_i}(x).
\end{equation*}
\end{proposition}

From the $VN$-spaces, one can construct the Voronoi polytopes:
\begin{theorem}
Given a polyhedral norm ${\mathcal L}$, suppose that we have a
tiling by $VN$-spaces $(P_i)_{i\in I}$.
Then:

(i) For all points $v\in \ZZ^n$, it holds
\begin{equation*}
V_{\leq}(N_{{\mathcal L}},v)=\cup \left\{P_i\mbox{~with~} v\in Near(P_i) \right\}.
\end{equation*}

(ii) For all points $v\in \ZZ^n$, it holds
\begin{equation*}
V_{ < }(N_{{\mathcal L}},v)=\cup \left\{P_i\mbox{~with~} \{v\} = Near(P_i) \right\}.
\end{equation*}

\end{theorem}
\proof This is clear from the definitions. \qed

Since the $VN$-spaces are polytopal, they also can be described by their vertices. On the other hand, because of a degree of arbitrariness in the choice of the $VN$-spaces, we need a notion of vertices that is independent of the chosen partition into $VN$-spaces.

\begin{definition}\label{DEF_D_point}
A point $x_0\in \RR^n$ is called a {\em $D$-point} if it satisfies one of the following equivalent conditions:

(i) $x_0$ is a local maximum of $d_{min}(\dot, {\mathcal L})$

(ii) For all $VN$-spaces $P$, containing $x_0$, $\alpha_P$
attains its maximum on $x_0$.
\end{definition}
The equivalence is clear. The notion of $D$-point is inspired
by Delaunay polytope.
A Delaunay polytope (cf. Section \ref{SEC_Intro}) is the center of an
empty sphere in classical Voronoi theory.
This center $c$ is then a local maximum for the function $d_{min}$.
However, unlike the case of Delaunay polytopes, $D$-points are not
necessarily isolated. This is apparent for the $L^{\infty}$ norm on $\ZZ^2$
for which the Voronoi polytope is $[-1/2,1/2]^2$ and every point on the
boundary is at distance $1/2$ from a point of $\ZZ^2$ and so is a $D$-point.

The set of all $D$-points is an union of distinct polytopes from each $VN$-space.
On the other hand, the dimension of the set of $D$-points is a useful
invariant. We cannot say anything a priori on the topology of this
point-set.

We now define the notion of {\em vertex} for the Voronoi polytopes that
we are considering.
\begin{definition}\label{DEF_True_Vertex}
Given a point $x\in V_{<}({\mathcal L}, v)$ and a $VN$-space decomposition
of the space $\RR^n$, we say that $x$ is a {\em vertex} if there exist a number
$r$ of $(n-1)$-dimensional polytopes $H_1$, \dots, $H_r$ such that:

(i) $x$ belongs to all $H_i$.

(ii) Any $(n-1)$-dimensional polytope $K_i\subset H_i$ with $x\in K_i$ is contained in a unique $VN$-space.

(iii) If $n_i$ is the normal vector to $H_i$, then the rank of $(n_1, \dots, n_r)$ is equal to $n$.
\end{definition}
The first and third condition means that $x$ is uniquely determined by the faces in which it is contained; it is the same condition as for polytopes.
The second condition means that $H_i$ are real hyperplanes in the sense 
that they are not hyperplanes separating two $VN$-spaces contained in
the same Voronoi polytope.
We have to use $K_i$ in order to deal with the fact that, possibly, the $VN$-spaces do not define a face-to-face tiling of $\RR^n$. As a consequence, this notion of vertex is independent of the chosen $VN$-space decomposition.

\begin{definition}
For a polyhedral metric ${\mathcal L}$ and the lattice $\ZZ^n$, the {\em point
group} $Pt(\ZZ^n, {\mathcal L})$ is the group of matrices $A\in \GL_n(\ZZ)$
such that for any $\ell \in {\mathcal L}$, the function $\ell_A$ with $\ell_A(x)= \ell(Ax)$
belongs to ${\mathcal L}$. The point group is always a finite group.

The {\em affine linear symmetry group} $\Aut(\ZZ^n, {\mathcal L})$ of
${\mathcal L}$ for the lattice $\ZZ^n$
is the group generated by the point group and the translations along $\ZZ^n$.
\end{definition}

In order to use $VN$-spaces in the enumeration, let us prove a number
of properties for them.

\begin{theorem}
Any $VN$-space is bounded.
\end{theorem}
\proof The function $d_{min}$ is bounded from above by the covering radius $\cov({\mathcal L})$.
However, the function $N_{{\mathcal L}}(x-v)$ is unbounded: so, any $VN$-space is bounded as well.
Hence, a given point $v$ can be contained in only a finite number of $VN$-spaces. \qed

The difficulty that one faces is that in the definition of the $VN$-spaces,
we have to account for every case.
The inequalities of the form $\ell'(x-v) \leq \ell(x-v)$
correspond to the function $\ell(x-v)$ defining the norm $N_{\mathcal L}(x-v)$.

However, inequalities of the form $N_{\mathcal L}(x-w) \geq \ell(x-v)$ are more
problematic.
We need to select a form $\ell_w$ such that $N_{\mathcal L}(x-w) = \ell_w(x-w)$.
This gives us two sets of inequalities:
\begin{enumerate}
\item $\ell'(x-w) \leq \ell_w(x-w)$ for $\ell'\in {\mathcal L}$
\item and $\ell_w(x-w)\geq \ell(x-v)$.
\end{enumerate}
These inequalities are quite complex. We can have the second kind of inequality
redundant for a $w_0$ but the first kind of inequalities for $w_0$ defining
facets of the $VN$-space.
On the other hand, if we remove the first kind of inequalities for $w_0$,
then the second kind of inequalities could be violated.

So, there is no simple way of choosing the set $w$ of inequalities that can
define a $VN$-space. 
But a finite set is sufficient. 
Also, the lack of a clear cut way of definition prevents
the construction of tilings and provable algorithms.

However, if one has a polytope $P$, then 
by using the algorithm of Subsection \ref{SEC_ClosestPoint}, one
can test efficiently whether or not $P$ is a $VN$-space.
This allows us to write a program that can build some $VN$-space
objects. But we cannot at this point guarantee that the programs will
return tilings and that if they form a tiling, it is face-to-face.

Our approach below is to take all possible inequalities defined by all
vectors. This allow us to build a procedure that works in the 
considered case of rational polyhedral metrics.

Let us first examine the geometric structure of the determining inequalities.
An {\em affine hyperplane arrangement} (AHA) in $\RR^n$ is a family of hyperplanes
that belongs to a finite number $p$ of translation classes.
In other words, there exist some hyperplanes $H_i$ and vectors $v_i$
for $1\leq i\leq p$ such that any hyperplane in the class is of
the form $H_i + j v_i$ for $j\in \ZZ$.
The connected components of an AHA are called {\em cells}.

\begin{lemma}
For a rational linear form $\ell\in {\mathcal L}$, the set of hyperplanes 
\begin{equation*}
H_v=\{x\in \RR^n\mbox{~:~}\ell(x-v)=0\}\mbox{~for~} v\in \ZZ^n
\end{equation*}
is of the form $H_0 + j w$ with $j\in \ZZ$ and $w\in \ZZ^n$.
\end{lemma}
\proof Let $e_1, \dots, e_n$ be the standard basis of $\RR^n$ and write $\ell(v_i)$ in the form $m_i / d$ with $n_i\in \ZZ$ and $d\in \NN$.
Write $h=gcd(m_1, \dots, m_n)$. By B\'ezout theorem, there exists a vector $w\in \ZZ^n$, such that $\ell(v) = h/d$.
For any $v\in \ZZ^n$, we have $\ell(v) = r (h/d)$ with $r\in \ZZ$.
So, any hyperplane $H_v$ is, actually, a translate of $H_0$ by $r w$. \qed

The above lemma will be used for special AHA defined below.

\begin{definition}\label{DEF_adapted_AHA}
An AHA is called {\em adapted} to a polyhedral 
metric ${\mathcal L}$ if for every cell $E$ and every vector $v\in \ZZ^n$,
there exists a $\ell\in {\mathcal L}$ such that, for every
$\ell'\in {\mathcal L}$, we have $\ell'(x-v)\leq \ell(x-v)$ for all $x\in E$.
\end{definition}

What we want is that each $VN$-space is contained into the cell
of an adapted AHA.
Of course, one has first to prove the 
existence of such arrangement.
It is also preferable to have simpler AHA
that are easier to work with computationally.

\begin{theorem}
(i) For any polyhedral metric ${\mathcal L}$, the set ${\mathcal L} - {\mathcal L}$ is adapted.

(ii) For any symmetric polyhedral metric ${\mathcal L}$, the set $\{{\mathcal L} - {\mathcal L}\}\diagdown \RR {\mathcal L}$ is adapted.

(iii) For the $L^1$ metric on $\RR^n$, the set $\{e_1^*, \dots, e_n^*\}$ is adapted.

(iv) For the $L^{\infty}$, metric on $\RR^n$ the set $\{\pm e_i^* \pm e_j^*\}_{1\leq i<j\leq n}$ is adapted.
\end{theorem}
\proof (i) If we take all the hyperplanes $\ell(x-v) = \ell'(x-v)$,
then for any cell of the corresponding arrangements 
we have either $\ell(x-v) < \ell'(x-v)$ or the reverse.
Hence, the $\ell(x-v)$ are ordered, and so, totally ordered and 
this ordering is independent of $x$. 
Hence, there exist a $\ell$ such that $\ell(x-v)$ dominates the other 
values. So, ${\mathcal L} - {\mathcal L}$ is adapted.

(ii) Let us write ${\mathcal L} = \{\pm \ell_1, \dots, \pm \ell_p\}$.
For us write ${\mathcal S} = \{{\mathcal L} - {\mathcal L}\}\diagdown 
\RR {\mathcal L}$ and take a cell $E$ of the corresponding AHA.
For any $1\leq i\leq p$, $x\in E$ and $v\in \ZZ^n$, write
$I_i(x,v)=[-\vert \ell_i(x-v)\vert, \vert \ell_i(x-v)\vert ]$.
The fixed inequalities between the $\ell_i$ ensures, that for any $i<j$ and $v\in \ZZ^n$ we have either $I_i(x,v)\subset I_j(x,v)$ or $I_j(x,v)\subset I_i(x,v)$ for all $x\in E$.
The intervals are totally ordered; so, there exists an $i_0$ such that $I_i(x,v)\subset I_{i_0}(x,v)$ for all $x\in E$ and $1\leq i\leq p$.
Let us take $i\not= i_0$. There exists $\epsilon\in \{1, -1\}$ such that $\pm \ell_i(x-v) < \epsilon \ell_{i_0}(x-v)$ for all $x\in E$.
By summing the two inequalities, we get $0< \epsilon \ell_i(x-v)$ for all 
$x\in E$. So, ${\mathcal S}$ is adapted. 

(iii) Let us take a cell $E$ of the AHA
determined by $\{e_1^*, \dots, e_n^*\}$.
Let us fix $v\in \ZZ^n$.
We have $e_i^*(x-v)$ of fixed sign $\epsilon_i$ over the cell $E$.
So, the inequality $\ell= \sum_{i=1}^n \epsilon_i e_i^*$ dominates all others.

(iv) Follows from (ii) and Equation \eqref{L_Norm_Linfini}. \qed

If we take the set ${\mathcal L} - {\mathcal L}$, then
on any cell $C$ of the corresponding AHA
the order of the values 
$\{\ell(x-v)\}_{\ell\in {\mathcal L}}$ does not depend only on $x\in P$.
This is, actually, more than what we require for the $VN$-spaces,
since for each vector we only need one $\ell$ such that
$\ell'(x-v)\leq \ell(x-v)$ for all $\ell'$.

The enumeration algorithm, that will be designed, will enumerate the 
$VN$-spaces corresponding to an adapted AHA.
Two such $VN$-spaces 
are called {\em adjacent} if their intersection is of dimension $n-1$.
The following is essential to the enumeration method:

\begin{theorem}
Let ${\mathcal L}$ be a polyhedral norm and ${\mathcal V}$ an adapted set of vectors for ${\mathcal L}$. Then there exist a family of $VN$-spaces $(P_i)_{i\in I}$ which form a face-to-face tiling of $\RR^n$ that finitely refines the tiling by the cells of the AHA defined by ${\mathcal V}$.
\end{theorem}
\proof Let us take a cell $E$ of the AHA
defined by ${\mathcal V}$.
Since $E$ is compact, there is a finite number of
points $v_1, \dots v_m\in \ZZ^n$ which are
at distance at most covering radius $\cov({\mathcal L})$
from any point of $E$.

By the definition of the adapted set ${\mathcal V}$,
 the functions $\phi_i(x) = N_{\mathcal L}(x-v_i)$ are linear on the 
cell $E$.
The tentative $VN$-spaces $P_i$ are thus defined as
\begin{equation*}
P_i = \left\{x\in E\mbox{~:~} \phi_j(x)\leq \phi_i(x) \mbox{~for~}j\not= i\right\}.
\end{equation*}
The ones that are full-dimensional, determine a finite $VN$-space tiling
of $E$ and so, a tiling of $\RR^n$. \qed

\section{Algorithms}\label{SEC_Algorithms}

In \cite{DSV_voronoi}, a complete set of algorithms is developed for
computing with Euclidean metrics on high--dimensional lattices.
Here we build similar methods for polyhedral metrics by using $VN$-spaces.

\subsection{Closest point}\label{SEC_ClosestPoint}

In the Euclidean case, the key ingredient in the algorithm is the solution
of the closest vector problem, that is, for a given $x\in \RR^n$, to find
all points $v\in \ZZ^n$ minimizing $\Vert x - v \Vert$. The solution to
this problem is given by the Fincke-Pohst algorithm \cite{FinckePohst}.

For a given polyhedral norm ${\mathcal L}$ and distance $d$, the set of $v\in \ZZ^n$
such that $N_{\mathcal L}(x - v) \leq d$ corresponds to the integral points
of the following polytope
\begin{equation*}
P_{\mathcal L}(x, d)= \{v\in \RR^n \mbox{~:~}\ell(x - v)\leq d \mbox{~for~}\ell\in {\mathcal L}\}.
\end{equation*}
Thus, the solution of the same problem for polyhedral norms, i.e. computing
$d_{min}(x, {\mathcal L})$, can be solved if one can determine integer points
in a polytope.

An efficient algorithm for solving such problems is provided by the
software {\tt zsolve} available via \cite{4ti2}. Note that, in order to
have a faster program, we first try to minimize the value of $d$ by finding
a point $v\in \RR^n$ which is near to $x$, though not necessarily the nearest,
by small coordinate changes.

Another algorithm for which {\tt zsolve} is useful is when we want to test that a given polytope $\SymbVN$ is a $VN$-space:

\begin{flushleft}
\smallskip
\textbf{Input:} a polyhedral metric ${\mathcal L}$ and a polytope $\SymbVN$\\
\textbf{Output:} If $\SymbVN$ is a $VN$-space return {\bf true} and otherwise a certificate that it is not.\\
\smallskip
$E\leftarrow$ set of vertices of $\SymbVN$.\\
$c\leftarrow$ isobarycenter of $E$.\\
$v_0\leftarrow$ nearest point to $c$.\\
$\ell_0\leftarrow$ the form $\ell\in {\mathcal L}$ realizing the maximum of $\ell(c - v_0)$.\\
\textbf{if} the inequalities $\ell(x- v_0) \leq \ell_0(x-v_0)$ are not valid on $\SymbVN$ \textbf{then}\\
\hspace{2ex} {\bf return} a $\ell$ and $x\in \SymbVN$ satisfying $\ell(x - v_0) > \ell_0(x-v_0)$.\\
\textbf{end if}\\
${\mathcal F}\leftarrow \emptyset$.\\
\textbf{for} $\ell \in {\mathcal L}$ \textbf{do}\\
\hspace{2ex} $h\leftarrow$ maximum of $\ell_0(x) -\ell(x)$ over $E$.\\
\hspace{2ex} ${\mathcal F}\leftarrow {\mathcal F}\cup \left\{\ell_0(v_0) \leq h + \ell(x)\right\}$.\\
\textbf{end for}\\
${\mathcal I}\leftarrow$ set of integral points of polytope defined by ${\mathcal F}$.\\
${\mathcal N}\leftarrow \emptyset$\\
${\mathcal F}\leftarrow$ facets of $P$.\\
\textbf{for} $v\in {\mathcal I}$ \textbf{do}\\
\hspace{2ex} ${\mathcal F}'\leftarrow {\mathcal F}$.\\
\hspace{2ex} \textbf{for} $l\in {\mathcal L}$ \textbf{do}\\
\hspace{2ex} \hspace{2ex} ${\mathcal F}'\leftarrow {\mathcal F}' \cup \{\ell(x-v)\leq \ell_0(x-v_0)\}$.\\
\hspace{2ex} \textbf{end for}\\
\hspace{2ex} \textbf{if} polytope determined by ${\mathcal F}'$ is non-empty and full-dim. \textbf{then}\\
\hspace{2ex} \hspace{2ex} ${\mathcal N}\leftarrow {\mathcal N} \cup \{v\}$.\\
\hspace{2ex} \textbf{end if}\\
\textbf{end for}\\
\textbf{if} ${\mathcal N} = \emptyset$ \textbf{then}\\
\hspace{2ex} \textbf{return} {\bf true}\\
\textbf{end if}\\
\textbf{return} ${\mathcal N}$
\end{flushleft}
The idea of this algorithm is that we take upper bound on possible values of 
$\ell_0(x) - \ell(x)$ which gives a potentially larger polytope.
Then, for each of the integral point obtained by {\tt zsolve}, we check
if the intersection is non-trivial.

\subsection{Group algorithms for $VN$-spaces}\label{ALGO_Group}

Given a polytope defined by linear inequalities, it is well known that one can 
obtain an interior point by using linear programming and so, the problem
can be solved in polynomial time.
However, this is insufficient for some polyhedral enumeration, since one
would like to get a point that is in fact canonical, i.e. invariant under
affine transformations.
That is, we need a function $f_{can}$ from the set of
polytopes in $\RR^n$ to $\RR^n$ such that, for any affine transformation
$\phi$ of $\RR^n$ and polytope $\SymbVN$, it holds
$\phi(f_{can}(\SymbVN)) = f_{can}(\phi(\SymbVN))$.
No general polynomial time solution of this problem is known.

But in our case, one can simply compute the vertices
of the considered $VN$-space $\SymbVN$ and then take their isobarycenter $Iso(\SymbVN)$.
This is, of course, relatively expensive, but reasonable for the cases
considered.

This isobarycenter can then be used to test equivalence of $VN$-spaces.
Two $VN$-spaces $\SymbVN$ and $\SymbVN'$ are equivalent
if and only if $Iso(\SymbVN)$ is equivalent to $Iso(\SymbVN')$.
The stabilizer of an $VN$-space $\SymbVN$ is found to be equal
to the stabilizer of $Iso(\SymbVN)$.
Hence, one can apply the algorithms developed in \cite{DSV_voronoi}, compute
stabilizers and test equivalence.

\subsection{Finding an initial $VN$-space}\label{SEC_Initial_VNspace}

We first give an algorithm that is fundamental to our enumeration methods.
It takes a point $x_0$ and returns the full-dimensional $VN$-space
$\SymbVN$ that contains $x_0$ in its interior if $\SymbVN$ exists.

\begin{flushleft}
\smallskip
\textbf{Input:} a polyhedral metric ${\mathcal L}$ and a point $x_0\in \RR^n$ and an AHA ${\mathcal H}$.\\
\textbf{Output:} $VN$-space $\SymbVN$ if $x_0$ is in the interior of an $VN$-space, fail otherwise.\\
\smallskip
${\mathcal S} \leftarrow$ a set of vectors that span $\ZZ^n$ and is antipodal invariant.\\
${\mathcal C}\leftarrow $ the set of points of $\ZZ^n$ closest to $x_0$.\\
\textbf{if} ${\mathcal C}$ has more than one element \textbf{then}\\
\hspace{2ex} return {\bf fail} \\
\textbf{else}\\
\hspace{2ex} call $v_0$ this element.\\
\textbf{end if}\\
\textbf{if} there are two $\ell \in {\mathcal L}$ realizing $\max_{\ell\in {\mathcal L}} \ell(x_0-v_0)$ \textbf{then}\\
\hspace{2ex} return {\bf fail} \\
\textbf{else}\\
\hspace{2ex} call $\ell_0$ this element.\\
\textbf{end if}\\
\textbf{while} \textbf{do}\\
\hspace{2ex} ${\mathcal F} \leftarrow \emptyset$\\
\hspace{2ex} \textbf{for} $\ell \in {\mathcal L}$ \textbf{do}\\
\hspace{2ex} \hspace{2ex} ${\mathcal F} \leftarrow {\mathcal F} \cup \{\ell(x-v_0)\leq \ell_0(x-v_0)\}$\\
\hspace{2ex} \textbf{end for}\\
\hspace{2ex} \textbf{for} $w \in {\mathcal S}$ \textbf{do}\\
\hspace{2ex} \hspace{2ex} \textbf{if} there are two $\ell\in {\mathcal L}$ realizing $\max_{\ell \in {\mathcal L}} \ell(x_0-w)$ \textbf{then}\\
\hspace{2ex} \hspace{2ex} \hspace{2ex} return {\bf fail} \\
\hspace{2ex} \hspace{2ex} \textbf{else}\\
\hspace{2ex} \hspace{2ex} \hspace{2ex} call $l'$ this element\\
\hspace{2ex} \hspace{2ex} \textbf{end if}\\
\hspace{2ex} \hspace{2ex} ${\mathcal F} \leftarrow {\mathcal F} \cup \{\ell'(x-w)\geq \ell_0(x-v_0)\}$\\
\hspace{2ex} \hspace{2ex} \textbf{for} $\ell\in {\mathcal L}$ \textbf{do}\\
\hspace{2ex} \hspace{2ex} \hspace{2ex} ${\mathcal F} \leftarrow {\mathcal F} \cup \{\ell(x-w)\leq \ell'(x-w)\}$\\
\hspace{2ex} \hspace{2ex} \textbf{end for}\\
\hspace{2ex} \textbf{end for}\\
\hspace{2ex} $\SymbVN \leftarrow$ the convex bodies defined by the inequalities of ${\mathcal F}$.\\
\hspace{2ex} \textbf{if} $\SymbVN$ is a bounded convex polytope \textbf{then}\\
\hspace{2ex} \hspace{2ex} ${\mathcal N}\leftarrow$ set of points $w\in \ZZ^n$ for which the set\\
\hspace{2ex} \hspace{2ex}  defined by $N_{\mathcal L}(x - w)\leq \ell_0(x-v_0)$ intersects $\SymbVN$ nontrivially.\\
\hspace{2ex} \hspace{2ex} \textbf{if} ${\mathcal N}\subset {\mathcal S}$ \textbf{then}\\
\hspace{2ex} \hspace{2ex} \hspace{2ex} \textbf{if} $\SymbVN$ is not split by any hyperplane in ${\mathcal H}$ \textbf{then}\\
\hspace{2ex} \hspace{2ex} \hspace{2ex} \hspace{2ex} return $\SymbVN$.\\
\hspace{2ex} \hspace{2ex} \hspace{2ex} \textbf{end if}\\
\hspace{2ex} \hspace{2ex} \textbf{end if}\\
\hspace{2ex} \textbf{end if}\\
\hspace{2ex} ${\mathcal S} \leftarrow {\mathcal S} + {\mathcal S}$\\
\textbf{end while}\\
\end{flushleft}

The method for finding an initial point is the following. Take a 
non-zero random vector $v\in \RR^n$ and 
divide $v$ by integers $k>0$ until the closest point to $v/k$ 
is $0$. Then we use the above 
algorithm to find the initial $VN$-space $\SymbVN$. If it 
fails, then we take another random vector $v$ and iterate.
In the last loop ${\mathcal S} + {\mathcal S}$ is a {\em Minkowski sum}, i.e.
we take all the sums $s + s'$ with $s,s'\in {\mathcal S}$.

\subsection{Finding adjacent $VN$-spaces}

We outline here our adjacent $VN$-space finding algorithm.
Like the preceding ones, it is based on an iterative scheme:

\begin{flushleft}
\smallskip
\textbf{Input:} a polyhedral metric ${\mathcal L}$, a $VN$-space $\SymbVN$ and a facet $F$ of $\SymbVN$.\\
\textbf{Output:} the $VN$-space $\SymbVN'$ adjacent to $\SymbVN$ on $F$.\\
\smallskip
$e\leftarrow$ isobarycenter of the vertices of $\SymbVN$ contained in $F$.\\
$v\leftarrow$ vector pointing from $e$ to the exterior of $\SymbVN$.\\
$\lambda\leftarrow 1$\\
\textbf{while} \textbf{do}\\
\hspace{2ex} $x\leftarrow e + \lambda v$\\
\hspace{2ex} $\SymbVN' \leftarrow$ result of algorithm of Section \ref{SEC_Initial_VNspace} for $x$.\\
\hspace{2ex} \textbf{if} $\SymbVN'$ is different from {\bf fail} and has $F$ as a facet \textbf{then}\\
\hspace{2ex} \hspace{2ex} return $\SymbVN'$\\
\hspace{2ex} \textbf{end if}\\
\hspace{2ex} $\lambda\leftarrow \lambda/2$\\
\textbf{end while}\\
\end{flushleft}

\subsection{The full enumeration algorithm}

Here we put the various pieces of the sub-algorithms together and get our main algorithm.
Its structure is similar to the enumeration algorithm used for Delaunay polytopes in \cite{DSV_voronoi}.
It is also a variant of the Voronoi algorithm for enumerating perfect forms \cite{EnumDim8}.

\begin{flushleft}
\smallskip
\textbf{Input:} Polyhedral metric ${\mathcal L}$.\\
\textbf{Output:} Set~${\mathcal R}$ of all inequivalent full-dimensional $VN$-spaces for ${\mathcal L}$.\\
\smallskip
$T \leftarrow $ initial $VN$-space for ${\mathcal L}$.\\
${\mathcal R} \leftarrow \emptyset$.\\
\textbf{while} there is a $\SymbVN \in T$ \textbf{do}\\
\hspace{2ex} ${\mathcal R} \leftarrow {\mathcal R} \cup \{ \SymbVN \}$.\\
\hspace{2ex} $T \leftarrow T \setminus \{\SymbVN \}$.\\
\hspace{2ex} ${\mathcal F} \leftarrow \mbox{facets of $\SymbVN$}$.\\
\hspace{2ex} \textbf{for} $F \in {\mathcal F}$ \textbf{do}\\
\hspace{2ex} \hspace{2ex} Find full-dimensional $VN$-space $\SymbVN'$ adjacent to $\SymbVN$ on $F$.\\
\hspace{2ex} \hspace{2ex} \textbf{if} $\SymbVN'$ is not equivalent to an $VN$-space in ${\mathcal R}\cup T$ \textbf{then}\\
\hspace{2ex} \hspace{2ex}\hspace{2ex} $T \leftarrow T \cup \{\SymbVN'\}$.\\
\hspace{2ex} \hspace{2ex} \textbf{end if}\\
\hspace{2ex} \textbf{end for}\\
\textbf{end while}\\
\end{flushleft}

The orbits of facets of $\SymbVN$ are computed with respect to the
stabilizer of $\SymbVN$ computed from Subsection \ref{ALGO_Group}.
The technique is to use the polyhedral enumeration program introduced
in \cite{EnumDim8}.

Two checks are available for the computation.
The first one: given a full-dimensional $VN$-space $\SymbVN$,
take a random point $x$ in the interior of $\SymbVN$ and 
compute the containing full-dimensional $VN$-space.
If it is distinct from $\SymbVN$, then there is an error.
Another check comes from the volume formula. Suppose that we have $m$ orbits of 
full-dimensional $VN$-spaces of representative $\SymbVN_1$, \dots, $\SymbVN_m$.
Denote by $\vert O_i\vert$ the number of translation classes.
Then we have the formula
\begin{equation*}
1 = \sum_{i=1}^m \vert O_i\vert \tvol(\SymbVN_i) \mbox{~with~} 
\vert O_i\vert = \frac{\vert G_{pt}\vert}{\vert \Stab(\SymbVN_i)\vert}.
\end{equation*}
Our algorithm can be adapted with minimal modifications
to more crystallographic applications of finding the Voronoi cells
for a polyhedral metric and a point set of the form $\{\ZZ^n + v_i\}_{1\leq i\leq N}$.
What is a priori more problematic is to consider the general case
of a non necessarily rational metric ${\mathcal L}$.

\subsection{Related computations}
The $D$-points (cf. Definition \ref{DEF_D_point}) can be determined
in the following way.
Let us assume that the tiling defined by the $VN$-spaces is face-to-face.
Given a $VN$-space $\SymbVN$, we compute all its vertices.
By testing equivalence of points using algorithms of
Subsection \ref{ALGO_Group}, one can determine the orbits of
vertices and, in addition, the list of $VN$-spaces in which they are contained.
If a vertex $v$ realizes the maximum of $\alpha_\SymbVN$ in all cells in
which it is contained, then it is a $D$-point.
For each $VN$-space $\SymbVN$, we take the list of their vertices that
are $D$-points and their convex hull define a polytope.
The collection of all such polytopes define the $D$-points.

For finding the vertices of $\SymbVN$, we again use the $VN$-spaces and assume face-to-face tilings. We enumerate all points of $V_{<}$ coming as vertices of $VN$-space, take the collection of all the hyperplanes and then we do counting. The exterior planes are the ones that appear only once.

Both methods can be extended to non face-to-face tilings. The idea is to refine the relevant faces into a tiling of several faces on which one can apply previous methods.

\section{Applications}\label{SEC_Applications}

Below are given
two distinct applications that illustrate nicely above methods. We take the root lattices $\mathsf{A}_n$ and $\mathsf{D}_n$ in their natural embedding in $\RR^{n+1}$ and $\RR^n$.
We use the $L^1$ and $L^{\infty}$ polyhedral norms on $\RR^n$ and compute the full-dimensional $VN$-spaces for both lattices. In practical terms, the limit to the computation is $n=6$ and comes from the use of {\tt zsolve}, which is the limiting factor.

The lattice $\mathsf{A}_n$ is defined as
\begin{equation*}
\mathsf{A}_n = \left\{ x\in \ZZ^{n+1}\mbox{~:~} \sum_{i=1}^{n+1} x_i =0\right\}
\end{equation*}
and its point group is isomorphic to $\ZZ_2\times \Sym(n+1)$.
The lattice $\mathsf{D}_n$ is defined as
\begin{equation*}
\mathsf{D}_n = \left\{ x\in \ZZ^{n}\mbox{~:~} \sum_{i=1}^{n} x_i \equiv 0 \pmod 2\right\}
\end{equation*}
and its point group has size $2^{n}\times n!$ for both the $L^1$ and $L^{\infty}$ norms.

We compute a $VN$-space decomposition for $\mathsf{D}_n$ and $\mathsf{A}_n$ for $n\leq 6$ and $L^1$ and $L^{\infty}$.
As a result, we are able to state the following conjecture:

\begin{conjecture}
For both, $\mathsf{D}_n$ and $\mathsf{A}_n$, and for both, $L^1$ and $L^{\infty}$, it holds:

(i) The strict Voronoi polytope $V_{<}(N_{\mathcal L})$ is equal to the interior of $V_{\leq }(N_{\mathcal L})$.

(ii) The Voronoi polytope $V_{\leq}(N_{\mathcal L})$ is equal to the Voronoi polytope $V_{\leq}(N_{eucl})$ with $N_{eucl}$ being the standard Euclidean norm.
\end{conjecture}
In other words, the above conjecture states that $Vor_{\leq}(N_{\mathcal L})$ is the convex hull of its vertices. The list of vertices is given in \cite[pp. 206-207]{DL}. It also seems possible that the conjecture is valid for any $L^p$ norm with $1\leq p\leq \infty$.

Similar results hold and are easy to prove for the lattice $\ZZ^n$ and the $L^p$ norms. The norm $N$ is then non-degenerate and the Voronoi body $V_{\leq}(N_{L^p}, 0)$ is then $[-1/2,1/2]^n$.

\end{document}